\documentclass[12pt,a4paper]{article}
\usepackage[ansinew]{inputenc}
\usepackage{latexsym,amsfonts,amsmath,amssymb,paralist}

\usepackage{url}

\newtheorem{theorem}{Theorem}[section]
\newtheorem{lemma}[theorem]{Lemma}
\newtheorem{cor}[theorem]{Corollary}
\newtheorem{con}[theorem]{Conjecture}
\newtheorem{pro}[theorem]{Proposition}

\newcommand{\pg}{{\rm PG}}
\newcommand{\gf}{{\rm GF}}
\newcommand{\ag}{{\rm AG}}
\newcommand{\q}{{\rm Q}}

\newcommand{\PP}{{\cal P}}

\newcommand{\oo}{{\cal O}}
\newcommand{\cS}{{\cal S}}

\newcommand\C{\mathcal{C}}
\newcommand\B{\mathcal{B}}
\newcommand\K{\mathcal{K}}
\newcommand\cH{\mathcal{H}}
\newcommand\I{\mathrel{\mathrm I}}

\newcommand\IS{\cS=(\PP,\B,\I)}

\title{Direction problems in affine spaces} 
\author{Jan De Beule\thanks{The author is a postdoctoral research fellow of the Research Foundation Flanders -- Belgium (FWO).} }

\date{}

\begin{document}

\maketitle
\begin{abstract}
This paper is a survey paper on old and recent results on direction problems in finite dimensional affine spaces over a finite field.
\end{abstract}

\section{Introduction}

Let $p$ be prime and $q=p^h$, $h \geq 1$. Then $\gf(q)$ denotes the finite field of order $q$. 
A vector space of dimension $n$ over $\gf(q)$ will be denoted as $V(n,q)$. The $n$-dimensional
affine space over the field $\gf(q)$, denoted as $\ag(n,\gf(q))$, shortly, $\ag(n,q)$, is an incidence geometry
of which the elements are the additive cosets of the vector subspaces of $V(n,q)$. The incidence is symmetrised 
set-theoretic containment. As such, a point of $\ag(n,q)$ is represented by a unique vector of $V(n,q)$, a line is 
represented by a unique coset of a $1$-dimensional vector subspace, etc.
If $S$ is an element of the affine space, so represented by a coset of a vector subspace $\pi$, then the (affine) dimension
of $S$ is defined as the vector dimension of $\pi$. Hence, the elements of $\ag(n,q)$ have dimension $i$, $0 \leq i \leq n-1$,
and we call the elements of $\ag(n,q)$ also {\em affine subspaces}.

The $n$-dimensional projective space over the field $\gf(q)$, denoted as $\pg(n,\gf(q))$, shorly, $\pg(n,q)$, is an incidence
geometry of which the elements are the vector subspaces of the $(n+1)$-dimensional vector space $V(n+1,q)$. 
Incidence is symmetrised containment. As such, a point of $\pg(n,q)$ is represented by a unique vector line, and hence up to
a $\gf(q)$ scalar, by a unique non-zero vector of $V(n+1,q)$. The zero vector represents the projective {\em empty subspace}.
We call the elements of $\pg(n,q)$ projective subspaces, and the projective dimension of a projective subpace $S$ is one less
than its vector dimension. 

Consider now an affine space $\ag(n,q)$. It is well known that in $\ag(n,q)$ parallelism exists: given a point $P$ and an 
$i$-dimensional subspace $S$ not incident with $P$, then there exists a unique $i$-dimensional subspace $S'$ containing $P$
and not meeting $S$. We call $S$ and $S'$ parallel. It simply means that $S$ and $S'$ are cosets of the same $i$-dimensional 
vector subspace $\pi$ of $V(n,q)$. The {\em direction at infinity} of $S$ (and $S'$) is represented by this vector subspace $\pi$. 
The set of directions of all affine spaces constitutes in fact a projective space $\pg(n-1,q)$. The projective points of such a direction
are simply called the points at infinity of the affine subspace $S$. Clearly, a line of $\ag(n,q)$ has a point at infinity. We will often
denote the projective space at infinity as $\pi_{\infty}$.

Let $U$ be a point set in $\ag(n,q)$. We call a point $t \in \pi_{\infty}$ a {\em direction determined by
$U$} if there exist two different points $r,s \in U$ such that the affine line determined by $r$ and $s$
contains at infinity the point $t$. It is clear that a point set $U$ will determine a set of points at infinity, we denote
the set of {\em directions determined by $U$} as the set $U_D$. In the literature, when dealing with point sets in $\ag(2,q)$, 
sometimes a determined direction is also called a {\em determined slope}, since a determined direction is the point at infinity of a line,
and this point at infinity can indeed be represented by the slope of the line, and all points at infinity can be represented by all possible
slopes (i.e. the elements of $\gf(q) \cup \{\infty\}$). The following research questions have been addressed.

\begin{enumerate}
\item What are the possible sizes of $U_D$ given that $|U|=q^{n-1}$? What is the possible structure of $U_D$?
\item What are the possible sets $U$, $|U|=q^{n-1}$, given that $U_D$ (or its complement in $\pi_{\infty}$) or only $|U_D|$ is known?
\item Given that a set $N$ of directions is not determined by a set $U$, $|U|=q^{n-1}-\epsilon$, can $U$ be extended to a set $U'$, 
$|U'|=q^{n-1}$, such that $U'$ does not determine the given set $N$?
\end{enumerate}

Directions problems have been studied in the past and also in affine spaces over arbitrary (non-finite) fields. A notorious
example is \cite{Besicovitch:1963ly}. A derived problem over finite fields has been addressed in \cite{Dvir:2009ve} and 
\cite{Blokhuis:2008qf}. We will focus on direction problems in affine spaces over finite fields. The earliest reference is probably the book of 
L. R\'edei (\cite{Redei:1970kx}). This brings us seamlessly to the most used technique to study direction problems: the so-called 
polynomial method. However, this  paper is mainly meant to be a survey paper on old and recent results on (and some applications of) 
direction problems in affine spaces. A recent and detailed survey paper on the polynomial method is e.g. \cite{Ball:2012}.

\section{Results in the affine plane $\ag(2,q)$}\label{sec:planar}

Probably the first result to be mentioned is the following theorem. It addresses research question 1, and 
it appeared in \cite{Redei:1970kx}, where it is stated in terms of the set of difference quotients of a function 
on a finite field, and it is obtained as a (non-trivial) consequence of the theory of so-called lacunary polynomials 
over finite fields. 

\begin{theorem}\label{th:redei1}
Let $U$ be a subset of the affine plane $\ag(2,p)$, $p$ prime, such that $|U|=p$ and $U$ is not the set of points of a line. Then $U$
determines at least $\frac{p+3}{2}$ directions.
\end{theorem}

An alternative proof, including a characterization of the equality case, can be found in \cite{Lovasz:1983uq}. 
This characterization can be formulated as follows. Note that the following theorem is addressing research 
question 2.

\begin{theorem}[\cite{Lovasz:1983uq}]\label{th:ls1983}
For every prime $p>2$, up to affine transformation there is a unique set of p points in $\ag(2,p)$ determining $\frac{p+3}{2}$ directions. Up to affine transformation,
the point set is the graph of the function $f(x): \gf(p) \rightarrow \gf(p): x \mapsto x^{\frac{p+1}{2}}$.
\end{theorem}

In some papers, Theorem~\ref{th:redei1} is referred to as a theorem of R\'edei and Megyesi. However, L\'aszl\'o Megyesi and L\'aszl\'o R\'edei 
have no joint paper.  In the first version of the manuscript of \cite{Redei:1970kx},  R\'edei had only the bound $\frac{p+1}{2}$, 
which was improved to $\frac{p+3}{2}$ after a discussion with Megyesi. Since the arguments for the improvement were in fact already 
present in \cite{Redei:1970kx}, Theorem~\ref{th:redei1} appears in the final version of the book,  and no follow up paper by R\'edei and 
Megyesi was written\footnote{I am grateful to Tam\'as Sz\H{o}nyi for this clarification}. 

The work of R\'edei contains examples of point sets of which the number of determined direction can (easily) be computed. 
These examples are also attributed to Megyesi. For an integer $d \mid q-1, 1 < d < q-1$, let $G=\gf(q)^*$ 
be the multiplicative subgroup of $\gf(q)$ of order $d$. Set 
\[ U = \{(x,0): x \in G\} \cup \{(0,x):x\not \in G\}\,,\] 
then $U$ determines $q+1-d$ directions. For $q$ prime, this is the unique example characterized in Theorem~\ref{th:ls1983}. 

The book of R\'edei \cite{Redei:1970kx} is a treatise on so-called lacunary polynomials over finite fields. 
Its main objective is to characterize polynomials in one variable over a finite field under certain 
assumptions on its degree and on the absence of terms of high degree. The direction result is obtained 
as a consequence of the developed theory. 

It seems that the first improvement on R\'edei's original theorem only appears in 1995 in a paper by Blokhuis, 
Brouwer and Sz\H{o}nyi, \cite{Blokhuis:1995aa}. In 1999 these results are again improved by Blokhuis, Ball, 
Brouwer, Storme and Sz\H{o}nyi, \cite{Blokhuis:1999kq}. In 2003, an unresolved case in the main theorem of  
\cite{Blokhuis:1999kq} is closed by Ball in \cite{Ball:2003aa}. We state here the main theorem of \cite{Ball:2003aa}. 
The interested reader can compare the different versions of this theorem throughout the papers 
\cite{Blokhuis:1995aa, Blokhuis:1999kq, Ball:2003aa}. Note that the following theorem is stated in terms of functions 
on $\gf(q)$. Clearly, a point set $U$ of points of $\ag(2,q)$, not determining all directions, corresponds, 
up to affine transformation, always to the graph of a function $f: \gf(q) \rightarrow \gf(q)$:
\[ U_f = \{(x,f(x)): x \in \gf(q)\}\,. \]

Also the converse is true: the graph of a function will not determine all directions. 

\begin{theorem}[\cite{Ball:2003aa}, Theorem 1.1]\label{th:ball2003:1}
Let $f$ be a function from $\gf(q)$ to $\gf(q)$, $q=p^h$ for some prime $p$, and let $N$ be the number of directions determined by $f$. 
Let $s=p^e$ be maximal such that any line with direction determined by $f$ that is incident with a point of the graph of $f$ is incident with a multiple 
of $s$ points of the graph of $f$. One of the following holds:
\begin{compactenum}[(i)]
\item $s=1$ and $\frac{q+3}{2} \leq N \leq q+1$;
\item $\gf(s)$ is a subfield of $\gf(q)$ and $\frac{q}{s}+1 \leq N \leq \frac{q-1}{s-1}$;
\item $s=q$ and $N=1$.
\end{compactenum}
Moreover, if $s > 2$ then the graph of $f$ is $\gf(s)$-linear, and all possibilities for $N$ can
be determined explicitly (in principle).
\end{theorem}

The results in \cite{Blokhuis:1995aa,Blokhuis:1999kq,Ball:2003aa} are based on further elaboration of 
the techniques used in \cite{Redei:1970kx}. Essentially, the so-called R\'edei-polynomial is associated to the 
affine point set. Its algebraic properties are derived form the geometric conditions and vice versa. 
These papers, together with other papers on blocking sets, which will be mentioned further in this 
survey, can be considered as the founding papers of the so-called polynomial method in finite geometry. 

The original theorem of R\'edei (Theorem~\ref{th:redei1}) and its improvements in 
\cite{Lovasz:1983uq,Blokhuis:1995aa, Blokhuis:1999kq, Ball:2003aa} were further
elaborated by G\'acs and Ball. The following theorem is due to G\'acs, and appeared in \cite{Gacs:2003fk}, 
which is the continuation of work started in \cite{Gacs:1999aa,Gacs:2000aa}. Note that \cite{Gacs:2003fk} was 
submitted right after publication of \cite{Gacs:1999aa,Gacs:2000aa}. We give the original statement of the theorem.
 
\begin{theorem}[\cite{Gacs:2003fk}, Theorem 1.3]\label{th:gacs2003}
For every prime $p$, besides lines and the example characterized by Lov\'asz and Schrijver 
(see Theorem~\ref{th:ls1983}), any set of p points in $\ag(2,p)$ determines at least $\left[2 \frac{p-1}{3}\right]+1$ 
directions.
\end{theorem}

Recall the examples of Megyesi. When $3 \mid q-1$,  the set $U$ determines exactly $2 \frac{p-1}{3}+2$ directions. 
Hence Theorem~\ref{th:gacs2003} is almost sharp.

In \cite{Szonyi:1999kq}, the theorem of R\'edei is reviewed, with the focus on applications. A classification of the example can also be found there.
Based on this result, there is a slightly different proof of Theorem~\ref{th:ls1983} in \cite{Gacs:2007fj}, and a generalization as follows.

\begin{theorem}[\cite{Gacs:2007fj}, Theorem 1.3]\label{th:gacs}
Let $U$ be a subset of the affine plane $\ag(2,q)$, $q=p^2$, $p$ prime, such that $|U_D| \geq \frac{q+3}{2}$ directions. Then either
$U$ is, up to affine transformation, equivalent to the graph of the function $f(x): \gf(p) \rightarrow \gf(p): x \mapsto x^{\frac{q+1}{2}}$,
or $U$ determines at least $\frac{q+p}{2} +1$ directions.
\end{theorem}

In \cite{Polverino:1999nr}, an example attaining the bound of Theorem~\ref{th:gacs} is reached, which makes this bound sharp. Note that the
bound Theorem~\ref{th:gacs} is weaker but similar to the bound in Theorem~\ref{th:gacs2003}

Up to now, all results except for Theorem~\ref{th:ls1983} and Theorem~\ref{th:ball2003:1} (ii) give only information 
on the number of determined directions. A lot of attention has also been paid on characterization results of examples. 
There is however no sharp line between such characterization examples and
examples only giving information on the number of determined directions. In many papers, both go together. 
We first mention some older results that provide characterizations. 

Assume that $p$ is prime and that $f$ is any function from $\gf(p) \rightarrow \gf(p)$. Let $M(f)$ be the number of 
elements $c \in \gf(p)$ such that $x \mapsto f(x) + cx$ is a permutation of $\gf(p)$, which is equivalent with saying the $c$ is
a non-determined direction of the graph of $f$. 

It should be noted that permutation polynomials have been studied for their own interest. Let $f(x)$ be a permutation
polynomial over $\gf(q)$. The question for permutation polynomials $f(x)$ over $\gf(q)$ the polynomial $f(x)+cx$ 
is a permutation polynomial for many values $c \in \gf(q)$ is studied in \cite{Evans:1992fr}.

A result of Sz\H{o}nyi characterizes point sets contained in the union of 
two lines under certain assumptions on the determined directions. The following result of Sz\H{o}nyi is a kind of 
generalization of Theorem~\ref{th:ls1983}.

\begin{theorem}[\cite{Szonyi:1991ys}]
If $M(f) \geq 2$, and the graph of $f$ is contained in the union of two lines, then after affine transformation, 
the graph of $f$ is equivalent to the example of Megyesi.
\end{theorem}

One of the most recent, if not the most recent paper on the direction problem in the plane is \cite{Ball:2009fk}. 
To state the results, we have to introduce the notation $I(f)$, which was also used in \cite{Gacs:2003fk}. 

Consider again a function $f: \gf(p) \rightarrow \gf(p)$. By interpolation, any function determines a 
polynomial of degree at most $p-1$ over $\gf(p)$, and conversely, every such polynomial determines a function. 
A function $f$ corresponds with a polynomial $g(X) \in \gf(p)[X]$ of a certain degree. Clearly, the function
$f(X)^i$ for any $i$, is represented by $g(X)^i \in \gf(p)[X]$. But since $x^p-x = 0$ for all $x \in \gf(p)$, we may reduce 
$g(X)^i$ modulo $X^p-X$ to obtain a polynomial representation of $f^i(X)$. We call the degree of $g(X)^i$ modulo $X^p-X$ 
the degree of $f^i(X)$. Then $I(f)$ is defined as follows:
\[
I(f)=\min \left\{ i+j: \sum_{x \in \gf(p)} x^j f(X)^i \neq 0 \right\}\,.
\]
In \cite{Ball:2009fk} it is explained why for all $n \leq I(f)$ implies that $f(X)^i$ has degree at most $p-2-n+i$. The following results are then found.

\begin{theorem}[\cite{Ball:2009fk}, Theorem 2.4]\label{th:bg2009:1}
If $I(f)> \frac{p - 1 - 2\epsilon}{t} + t - 2 + \epsilon$ for some integer $t$ then every line meets the graph of f in at 
least $I(f)+3-t >\frac{p-1}{t} +1$ points or at most $t - 1$ points.
\end{theorem}

The authors state the following conjecture, based on the proof of the previous theorem.

\begin{con}\label{con:bg2009}
If $I(f)> \frac{p - 1 - 2\epsilon}{t} + t - 2 + \epsilon$ for some integer $t$, then the graph of $f$ is contained in an algebraic
curve of degree $t-1$.
\end{con}

The next theorem is actually the proof of this conjecture under extra assumptions.
\begin{theorem}[\cite{Ball:2009fk}, Theorem 2.6]\label{th:bg2009:2}
If $I(f)> \frac{p - 1 - 2\epsilon}{t} + t - 2 + \epsilon$ and there are $t-1$ lines incident with at least $t$
points of the graph of f then the graph of $f$ is contained in the union of these $t - 1$ lines.
\end{theorem}

Putting $t=2$ in Theorem~\ref{th:bg2009:1} yields the following corollary, which is a reformulation of Theorem~\ref{th:ls1983}.

\begin{cor}[\cite{Ball:2009fk}, Theorem 3.1]\label{co:bg2009:1}
If $I(f) \geq \frac{p+1}{2}$ then $f$ is linear.
\end{cor}

The next theorem is a generalization of a theorem in \cite{Gacs:2003fk}. 

\begin{theorem}[\cite{Ball:2009fk}, Theorem 3.2]\label{th:bg2009:3}
If $I(f) \geq \frac{p+5}{3}$ then the graph of $f$ is contained in an algebraic curve of degree $2$.
\end{theorem}

In \cite{Wan:1995zr}, information on $M(f)$ in terms of the degree of $f$ is obtained.

Up to now, all results where related to research questions 1 and/or 2. The following result is a stability result 
deals which addresses research question 3.

\begin{theorem}[\cite{Szonyi:1996fk}, Theorem 4]\label{th:szt}
A set $U$ of $q-k > q-\frac{\sqrt{q}}{2}$ points of $\ag(2,q)$ for which $|U_D| \leq \frac{q+1}{2}$, can be 
extended to a set $U'$ of $q$ points of $\ag(2,q)$ such that $U_D = U'_D$.
\end{theorem}

The case $q=p$ was handled separately in \cite{Szonyi:1999kq}, using lacunary polynomials. The proof of 
Theorem~\ref{th:szt} is much more dependent on algebraic geometric arguments. Several remarks with 
refinements and  consequences under particular assumptions are found in \cite{Szonyi:1996fk}. 
\begin{compactenum}[\rm (i)]
\item For $q$ prime, the bound $q-k > q-\frac{\sqrt{q}}{2}$ can be improved to $q-k > q- \frac{p+45}{20}$.
\item The case $k=1$ has a very short proof. Since $n=1$, each $y \not \in U_D$, there is a unique line $L_y$ not meeting
$U$. From the short argument of the proof, it is deduced that all lines $L_y$, $y \not \in U_D$, pass through a common point.
\end{compactenum}

The following theorem is also found in \cite{Szonyi:1999kq}. 

\begin{theorem}[\cite{Szonyi:1999kq}]
Let $U$ be a set of $k$ points of $\ag(2,p)$, $p$ prime, such that not all $k$ points of $U$ are collinear. 
Then $|U_D| \geq \frac{k+3}{2}$. 
\end{theorem}

Consider the affine plane $\ag(2,p)$, and consider a coset of a multiplicative subgroup $H \leqslant \gf(p)^*$. 
Set 
\[ U = \{(x,0): x \in H\} \cup \{(0,x):x \in H\} \cup \{(0,0)\}\,,\] 
then $|U|=2|H|+1$ and $|U_D| = |H|+2 = \frac{k+3}{2}$, hence, when $k = 2d+1$ with $d \mid p-1$, the 
bound in the theorem is sharp. This ``Megyesi-type'' example is due to Aart Blokhuis.

The most recent result on planar direction problems is found in \cite{Fancsali:2013oq}. This paper 
actually addresses a variation on research question 1 in the plane. The authors consider 
a set $U$ in $\ag(2,q)$ of less than $q$ points, and derive a result similar to the results in \cite{Blokhuis:1999kq} 
(which are part of Theorem~\ref{th:ball2003:1}).

Let $q=p^h$, $p$ prime. Let $U$ be a set of points of $\ag(2,q)$. Let $d$ be a direction at infinity, then define $s(d)$ as the greatest power
of $p$ such that each line $l$ of direction $d$ meets $U$ in zero modulo $s(d)$ points, and define
\[ s = \mathrm{min} \{s(y): y \in U_D\}\,. \]

Let $U = \{(a_i,b_i): 1 \leq i \leq |U|\}$. The R\'edei-polynomial associated to the set $U$ is defined as
\[ R(X,Y) = \prod_{i=1}^{|U|} (X-a_iY+b) = X^n + \sum_{j=0}^{n-1}\sigma_{n-j}(Y) X^j \]

The following proposition is in principle straightforward (as \cite{Fancsali:2013oq} is 
self contained, a proof can be found there). We describe it to introduce one more notion that 
is needed to formulate the main result from \cite{Fancsali:2013oq}).

\begin{lemma}
If $y \in U_D$, then $R(X,y) \in \gf(q)[X^{s(y)}] \setminus \gf(q)[X^{p\cdot s(y)}]$. If $y \not \in U_D$, then
$R(X,y) \mid X^q-X$. 
\end{lemma}

The above observation leads to the definition of a polynomial $H(X,Y)$ as follows. Consider $R(X,Y)$ as a univariate
polynomial over the ring $\gf(q)[Y]$. Since $R(X,Y)$ is monic, division with remainder of $X^q-X$ by $R(X,Y)$ yields
the quotient $Q(X,Y)$ and remainder $S(X,Y)$, define $H(X,Y) := - S(X,Y) - X$.  Properties of the polynomial
$H(X,Y)$ are then shown in  \cite{Fancsali:2013oq}, e.g. that $H(X,Y)$ is a constant polynomial if $|U_D| = 1$. Suppose 
that $|U_D| > 1$, then define $t(d)$ as the maximal power of $p$ such that $H5X,d) = f_d(X)^{t(d)}$ for some
$f_d(X) \not \in \gf(q)[X^p]$, and define

\[ t = \mathrm{min} \{s(y): y \in U_D\}\,. \]

The main theorem can now be formulated.

\begin{theorem}[\cite{Fancsali:2013oq}, Theorem 17]\label{th:fancsali13}
Let $U$ be an arbitrary set of points of $\ag(2,q)$. Assume that $\infty \in U_D$, then one of the following holds.
\begin{compactenum}[(i)]
\item $1=s \leq t < q$ and $\frac{|U|-1}{t+1} + 2 \leq |U_D| \leq q+1$, or,
\item $1 < s \leq t < q$ and $\frac{|U|-1}{t+1} + 2 \leq |U_D|  \leq \frac{|U|-1}{s-1} \leq q+1$, or,
\item $1 \leq s \leq t=q$ and $U_D = \{\infty\}$. 
\end{compactenum}
\end{theorem}

\section{Planar direction problems, blocking sets and the polynomial method}

A {\em blocking set} of a projective plane $\Pi$ is set of points $B$ such that any line of $\Pi$ meets $B$ in at least one point.
We call a blocking set {\em trivial} if it contains a line, and minimal if no point of $B$ can be deleted. The study of blocking sets
of Desarguesian projective planes in particular is important in finite geometry, and results on blocking sets of the Desarguesian
projective plane $\pg(2,q)$ have many applications in the study of other substructures in finite projective spaces and finite classical
polar spaces. We shortly describe in this section the connection between blocking sets, direction problems and the polynomial method.

Consider now a point set $U$ of size $q$ in the affine plane $\ag(2,q)$. The extension of $\ag(2,q)$ to the projective plane, by adding the
slopes at infinity as points is well known. If $U$ is the set of points of a line, then $U$ determines one directions $d$ $U \cup \{d\}$ is then
a projective line, indeed meeting all lines of $\pg(2,q)$. So assume that $U$ is not contained in an affine line, then the set $B := U \cup U_D$ 
is a minimal blocking set of $\pg(2,q)$. Consider any line $l$ of $\pg(2,q)$, not the line at infinity. If $l$ has a slope $d \in U_D$, then $l$ 
meets $U$ in at least two points. Furthermore, since $|U|=q$, there exists at least one line $m$ on $d$ not meeting $U$. Hence $d$ cannot
be removed from $B$. If $l$ has slope $d \not \in U_D$, then every line on $d$ meets $U$ in exactly one point. Since $|U| = q$ and there 
are $q+1$ lines on a point, on every point $P \in U$ there are at least two lines on $P$ meeting $U$ only in $P$. Finally, the line at infinity meets
$U_D$. Hence, $B$ is a minimal blocking set of size $q+n$, where $n := |U_D|$, and there exists at least one line meeting $B$ in exactly $n$ points.
We conclude that any point set $U$ of size $q$ gives rise to such a blocking set, which is called a {\em blocking set of R\'edei type.} 

The converse is not true: not every minimal blocking set is of R\'edei type, and so is not constructed as $U \cup U_D$. However, the idea of
R\'edei blocking sets, the associated direction problem, and the results of R\'edei have been inspiring to study small blocking sets. We first mention the
following result of Bruen.

\begin{theorem}[\cite{Bruen:1971zl}]\label{th:bruen}
Let $B$ be a blocking set of a finite projective plane $\Pi$ of order $n$. Then $|B| \geq n+\sqrt{n}+1$, and $|B| = n+\sqrt{n}+1$ if and only if $B$ 
is a Baer sub plane of $\Pi$.
\end{theorem}

An alternative proof, based on elementary counting techniques simplifies the proof Theorem~\ref{th:bruen}. It can be found in \cite{Bruen:1977aa}.

Considering the examples of Megyesi, and Theorem~\ref{th:ls1983}, it is clear that for $p$ odd prime a minimal blocking set of size 
$p+\frac{p+3}{2} = \frac{3(p+1)}{2}$ exists. It was only shown by Blokhuis in 1994 that no blocking sets of size smaller 
than $\frac{3(p+1)}{2}$ exists. 

\begin{theorem}[\cite{Blokhuis:1994aa}]
Assume that $q=p$ is an odd prime. If $B$ is a blocking set of $\pg(2,q)$ and $|B| \leq \frac{3p+1}{2}$, then $B$ contains all the points of a line
\end{theorem}

The proof of this theorem in \cite{Blokhuis:1994aa} is based on a generalization of a result of R\'edei on lacunary polynomials.  
The result was generalized to planes of prime power order in the following theorem.

\begin{theorem}[\cite{Blokhuis:1996ve}, Theorem 6]\label{th:blokje}
Assume that $q=p^{2e+1}$, $p$ and odd prime, $e \geq 1$. Then a minimal blocking set of $\pg(2,q)$ has size at least $q+\sqrt{pq}+1$ points.
\end{theorem}

The following examples of blocking sets are found in \cite{Redei:1970kx}.
Let $q=p^e$, $e > 1$, and let $\gf(q_1)$ be a subfield of $\gf(q)$. Using $f=T$, the trace function from $\gf(q) \rightarrow \gf(q_1)$, 
the graph of the function $f$ determines $\frac{q}{q_1}+1$ directions, so this construction yields a blocking set of size 
$q+\frac{q}{q_1}+1$. Hence, for $q=p^3$, the bound of Theorem~\ref{th:blokje} is sharp.

More information on blocking sets of Desarguesian projective 
planes can be found in \cite{BSS:2012}. The situation for non-Desarguesian planes is more complicated, especially 
the construction of examples. Apart from \cite{Bruen:1971zl} and  \cite{Bruen:1977aa}, the papers 
\cite{Bierbrauer:1980hl,Bierbrauer:1981qd,Kitto:1989qr} are interesting references for blocking sets of
non Desarguesian planes.

\section{Direction problems in affine spaces}

An early result on directions problems in affine spaces is found in \cite{Storme:2001nx}. It is motivated by the 
study of R\'edei type blocking sets in projective spaces. The main theorem is the following. Its proof uses
almost exclusively geometric and combinatorial arguments.

\begin{theorem}[\cite{Storme:2001nx}, Theorem 16]
Let $U$ be a set of points of $\ag(n,q)$, $n \geq 3$, $q=p^h$, $|U|=q^k$. Suppose that $U$ is a $\gf(p)$-linear set of points
and that $|U_D| \leq \frac{q+3}{2}q^{k-1} + q^{k-2} + \dots + q^2+q$. If $(n-1 \geq (n-k)h$, then $U$ is a cone with an
$(n-1-h(n-k))$-dimensional vertex at infinity and with a $\gf(p)$-linear point set $U_{(n-k)h}$ of size $q^{(n-k)(h-1)}$, contained
in some affine $(n-k)h$-dimensional subspace of $\ag(n,q)$.
\end{theorem}

The following theorem is found in Sziklai, \cite{Sziklai:2006eu}. Note that Theorem~\ref{th:gacs2003} is used 
to prove it.

\begin{theorem}[\cite{Sziklai:2006eu}, Theorem 15]
Let $U$ be a set of $p^2$ points in $\ag(2,p)$, $p$ prime, such that $|U_D| < \frac{2p(p-1)}{3} +2p$. 
Then the set $U$ is either a plane or a cylinder with the projective triangle as base if $p>11$.
\end{theorem}

In Section~\ref{sec:planar} we have seen that many direction results in affine planes are obtained by studying 
the graph of a  function in one variable over the finite field $\gf(q)$. The following result, which is probably (one of) 
the first results on direction problems in three dimensional affine spaces, generalizes this approach
to functions in two variables, and can be found in \cite{Ball:2006mz}. It addresses research question 2.

\begin{theorem}[\cite{Ball:2006mz}, Theorem 2.4]\label{th:ag3q1}
Let $q=p^h$, $prime$. If $U$ is a set of $q^2$ points of $\ag(3,q)$ that does not determine at least 
$p^eq$ directions for some $e \in \mathbb{N} \cup \{0\}$, then every plane meets $U$ in $0 \mod p^{e+1}$ points.
\end{theorem}

As in e.g.  \cite{Blokhuis:1995aa,Blokhuis:1999kq,Ball:2003aa}, a R\'edei-polynomial is associated to 
the point set in $\ag(3,q)$ determined as the graph of a function $f$ in two variables over $\gf(q)$. 
The R\'edei-polynomial is now a polynomial in three variables, again lacunary, from which again 
strong algebraic properties can be derived. The use of R\'edei-polynomials in more variables is 
further described in \cite{Ball:2004kq}.

Theorem~\ref{th:ag3q1} is generalized and improved in \cite{Ball:2008rc}. 

\begin{theorem}[\cite{Ball:2008rc}, Theorem 1.3]\label{th:ball2008}
Let $q=p^h$ and $1 \leq p^e < q^{k-2}$, where $e$ is a non-negative integer. If there are more than 
$p^e(q-1)$ directions not determined by a set $U$ of $q^{k-1}$ points in $\ag(k,q)$, then every hyperplane 
meets $U$ in $0 \mod p^{e+1}$ points.
\end{theorem}

The generalization of Theorem~\ref{th:ag3q1} to general dimension is relatively straightforward from the 
arguments used to show the theorem. The improvement is based on using the representation 
$\gf(q) \times \gf(q)^{k-1}$ for $\ag(k,q)$, and then associating a R\'edei-polynomial
to the set $U$, which will, due to the used representation, be a polynomial in only two variables. 

To construct a set of $q^2$ points in $\ag(3,q)$, that do not determine {\em many} directions, one could consider a 
set $U$ of $q$ points of $\ag(2,q)$, determining few directions, and then form a cylinder from $U$. The easiest example
is to take for $U$ the set of $q$ points on a line, then the corresponding set in $\ag(3,q)$ will be an affine plane. In \cite{Ball:2008rc},
it is noted that this procedure is actually the only known way of constructing such sets in $\ag(3,q)$. Connecting this with
Theorems~\ref{th:ag3q1}~and~\ref{th:ball2008}, the following conjectures are explained in \cite{Ball:2008rc}.

\begin{con}[\cite{Ball:2008rc}, Conjecture 5.1]
Let $U$ be a set of $p^2$ points in $\ag(3,p)$, $p$ prime, and let $N$ be the set of non-determined directions. If $|N| \geq p$,
then $U$ is the union of $p$ parallel lines.
\end{con}

By Theorem~\ref{th:ball2008}, it follows that a set satisfying the conditions of the conjecture, meets every plane
of $\ag(3,p)$ in $0 \mod p$ points. As such, the following conjecture is introduced in \cite{Ball:2008rc} as the {\em strong
cylinder conjecture}.

\begin{con}[\cite{Ball:2008rc}, Conjecture 5.2]\label{con:cylinder}
Let $U$ be a set of $p^2$ points in $\ag(3,p)$, $p$ prime. If $U$ has the property that every planes meets $U$
in $0 \mod p$ points, then $U$ is the union of $p$ parallel lines.
\end{con}

As a final note, a generalization of this conjecture in terms of finite abelian groups is given in \cite{Ball:2008rc}. 
In \cite{Ball:2012}, a weaker form of Conjecture~\ref{con:cylinder}, assuming on top that at least $p$ directions 
are not determined by the set $U$, is proposed as an open problem. 

The theorems and corollary mentioned in Section~\ref{sec:planar} (Theorems~\ref{th:bg2009:1}, 
\ref{th:bg2009:2}, \ref{th:bg2009:3} and Corollary~\ref{co:bg2009:1}) from \cite{Ball:2009fk} are purely planar 
results and continue on work started in \cite{Gacs:2003fk}. The work in \cite{Ball:2008uq}, contains some 
slight improvements of \cite{Gacs:2003fk}. Its planar results were then improved further 
in \cite{Ball:2009fk} (and are therefore not explicitly mentioned in Section~\ref{sec:planar}.
In this section, we discuss that a generalization of the planar results to three dimensional affine spaces 
in \cite{Ball:2008uq}.

Let $U$ be a set of $q$ points in $\ag(3,q)$. A line $l$ at infinity is called {\em not determined by $U$} if 
every affine plane through $l$ contains exactly one point of $U$. A point set in $\ag(3,q)$ will now be 
studied as the graph of a pair of functions over the finite field $\gf(q)$. Then a line not determined by the 
graph $\{(x,f(x),g(x)): x \in \gf(q) \}$, is corresponds with a pair $(c,d)$ such that $f(x) + cg(x) + dx$ is a permutation
polynomial.

Let $p$ be prime and let $f$ and $g$ be two functions over $\gf(p)$. Define
$M(f,g)$ be the number of pairs  $(a,b) \in \gf(p)^2$ such that $f(x) + ag(x) + bx$ is a permutation polynomial, and let
\[
I(f,g)=\min \left\{ k+l+m: \sum_{x \in \gf(p)} x^k f(X)^l g(x)^m \neq 0 \right\}\,.
\]

The main result of the non-planar part of \cite{Ball:2008uq} is then the following theorem.

\begin{theorem}[\cite{Ball:2008uq}, Theorem 3.2]\label{th:dbgacs}
Let $s = \lceil \frac{p-1}{6} \rceil$. If $M(f,g) > (2s+1)(p+2s)/2$, then there are elements $c,d,e \in \gf(p)$ 
such that $f(x)+cg(x)+dx+e = 0$.
\end{theorem}

As a final remark, an example is provided in \cite{Ball:2008uq} that shows that the bound on $M(f,g)$ is the right order of magnitude.

In \cite{Sziklai:2012dn}, the notion of direction at infinity is generalized in a different way to arbitrary subspaces at infinity. 
Let $U$ be a set of points of $\ag(n,q)$, and let $k \leq n-2$ be a fixed integer. A projective subspace $S$ of dimension 
$k$ at infinity is determined by $U$ if there is an affine subspace $T$ of dimension $k+1$ meeting the hyperplane at infinity 
in $S$, is spanned by the point set $U \cap T$. One observes easily that $|U| \leq q^{n-1}$ if not al projective subspaces of 
dimension $k$ at infinity are determined. The following theorem gives some information in the three-dimensional case.

\begin{theorem}[\cite{Sziklai:2012dn}, Theorem 7] Let $U$ be a set of $q^2$ points of $\ag(3,q)$. Let $U_L$ be the set of lines at 
infinity determined by $U$, and let $N$ be the set of non-determined lines at infinity. Then one of the following holds.
\begin{compactenum}[(i)]
\item The set $U$ determines all lines at infinity (so $|N|=0$);
\item $|N|=1$ and there is a parallel class of affine planes such that $U$ contains one (arbitrary) complete line in each of its planes;
\item $|N|=2$ and the set $U$ together with two undetermined lines at infinity form a hyperbolic quadric or $U$ contains $q$ parallel lines;
\item $|N| \geq 3$ and then $U$ contains $q$ parallel lines.
\end{compactenum}
\end{theorem}

Up to here, this section has been devoted only to results addressing research questions 1 and/or 2 (or variations on it). 
The following theorems deals with research question 3. It was first shown in \cite{De-Beule:2008fk} in three dimensions.
The following formulation is taken from \cite{Ball:2012}, where a proof for general $n$ is given using the representation
$\gf(q) \times \gf(q)^{n-1}$ for $\ag(n,q)$.

\begin{theorem}[\cite{Ball:2012}, Theorem 6.8]\label{th:dbgacs}
Let $q=p^h$, $p$ and odd prime. A set $U$ of $q^{n-1}-2$ points of $\ag(n,q)$ that does not determine a set $D$ of at 
least $p+2$ directions, can be extended to a set $U'$ of $q^{n-1}$ points, not determining the same set $D$ of directions.
\end{theorem}

The main motivation of the work in \cite{De-Beule:2013uq} was to study the problem of Theorem~\ref{th:dbgacs}
in an alternative way.

\begin{theorem}[\cite{De-Beule:2013uq}, Theorem 12]
Let $q=p^h$, $p$ prime. Let $U$ be a set of $q^2-\varepsilon$ points of $\ag(3,q)$, where $\varepsilon < p$. 
Put $N=\pi_\infty\setminus U_D$ the set of non-determined directions. Then $N$ is contained in a plane 
curve of degree $\varepsilon^4-2\varepsilon^3+\varepsilon$ or $U$ can be extended 
to a set $U'$, $|U'|=q^2$, $U_D = U'_D$.
\end{theorem}

\begin{theorem}[\cite{De-Beule:2013uq}, Theorem 13]\label{th:dbszt2013}
Let $n \geq 3$. Let $U\subset\ag(n,q)\subset\pg(n,q)$, $|U|=q^{n-1}-2$. Put $N=\pi_\infty\setminus U_D$ 
the set of non-determined directions. Then $U$ can be extended to a set $U'$, $|U'|=q^{n-1}$, 
$U_D = U'_D$, or the points of $N$ are collinear and $|N| \leq \lfloor \frac{q+3}{2} \rfloor$,
or the points of $N$ are on a (planar) conic curve.
\end{theorem}

\section{Direction problems in affine spaces and special point sets in generalized quadrangles}

A (finite) \emph{generalized quadrangle} (GQ) is an incidence structure
$\IS$ in which $\PP$ and $\B$ are disjoint non-empty sets of objects called
points and lines (respectively), and for which $\I \subseteq (\PP \times
\B) \cup (\B \times \PP)$ is a symmetric point-line incidence relation
satisfying the following axioms:
\begin{itemize}
\item[(i)] each point is incident with $1+t$ lines $(t \geqslant 1)$ and
two distinct points are incident with at most one line;
\item[(ii)] each line is incident with $1+s$ points $(s \geqslant 1)$ and
two distinct lines are incident with at most one point;
\item[(iii)] if $x$ is a point and $L$ is a line not incident with $x$,
then there is a unique pair $(y,M) \in \PP \times \B$ for which $x \I M \I
y \I L$.
\end{itemize}
The integers $s$ and $t$ are the parameters of the GQ and $\cS$ is said to
have order $(s,t)$. If $s=t$, then $\cS$ is said to have order $s$. If $\cS$
has order $(s,t)$, then $|\PP| = (s+1)(st+1)$ and $|\B| =
(t+1)(st+1)$ (see e.g. \cite{Payne:1984kx}). The \emph{dual} $\cS^D$ of a GQ
$\IS$ is the incidence structure $(\B,\PP,\I)$. It is again a GQ.

An {\em ovoid} of a GQ $\cS$ is a set $\oo$ of points of $\cS$ such that every line
is incident with exactly one point of the ovoid. An ovoid of a GQ of order
$(s,t)$ has necessarily size $1+st$. An {\em partial ovoid} of a
GQ is a set $\K$ of points such that every line contains {\em at most} one point
of $\K$. A partial ovoid $\K$ is called {\em maximal} if and only if $\K \cup
\{p\}$ is not a partial ovoid for any point $p \in \PP \setminus \K$, in other words, if
$\K$ cannot be extended. It is clear that any partial ovoid of a GQ of order
$(s,t)$ contains $1+st-\rho$ points, $\rho \geq 0$, with $\rho = 0$ if and only
if $\K$ is an ovoid. 

Consider a non-singular quadratic form $f$ acting on $V(5,q)$, which is, up to coordinate
transformation, unique. The points of $\pg(4,q)$ that are totally singular with relation to $f$ constitute 
the parabolic quadric $\q(4,q)$. Since $f$ has Witt index two, There exist projective lines on $\pg(4,q)$ that
are completely contained in $\q(4,q)$. It is well known that $\q(4,q)$ is actually a GQ of order $q$. This 
GQ is one of the so-called finite classical generalized quadrangles. The motivation of \cite{De-Beule:2008fk}
was to study the extendability of partial ovoids of $\q(4,q)$, $q$ odd. Using an alternative representation
of $\q(4,q)$, this extendability problem translates directly to a stability question on sets of size $q^2-\epsilon$ of
$\ag(3,q)$, not determining a given set of directions. 

An {\em oval} of $\pg(2,q)$ is a set of $q+1$ points $\C$, such that no three points
of $\C$ are collinear. When $q$ is odd, it is known that all ovals of 
$\pg(2,q)$ are conics. When $q$ is even, several other examples and infinite families are
known, see e.g. \cite{DCD:2012}. The GQ $T_2(\C)$ is defined as follows. Let
$\C$ be an oval of $\pg(2,q)$, embed $\pg(2,q)$ as a plane in
$\pg(3,q)$ and denote this plane by $\pi_{\infty}$. Points are defined as
follows:
\medskip
\begin{compactenum}[{\rm (i)}]
\item[(i)] the points of $\pg(3,q) \setminus \pg(2,q)$;
\item[(ii)] the planes $\pi$ of $\pg(3,q)$ for which $|\pi \cap \C| =
1$; 
\item[(iii)] one new symbol $(\infty)$.
\end{compactenum}
\medskip
Lines are defined as follows:
\medskip
\begin{compactenum}[{\rm (a)}]
\item[(a)] the lines of $\pg(3,q)$ which are not contained in $\pg(2,q)$
and meet $\C$ (necessarily in a unique point);
\item[(b)] the points of $\C$.
\end{compactenum}
\medskip
Incidence between points of type (i) and (ii) and lines of type (a) and
(b) is the inherited incidence of $\pg(3,q)$. In addition, the point
$(\infty)$ is incident with no line of type (a) and with all lines of type
(b). It is straightforward to show that this incidence structure is a
GQ of order $q$. The following theorem (see e.g. \cite{Payne:1984kx}) allows us to
use this representation. 

\begin{pro}
The GQs $T_2(\C)$ and $\q(4,q)$ are isomorphic if and only if $\C$ is a conic of
the plane $\pg(2,q)$.
\end{pro}

Suppose now that $\oo$ is a (partial) ovoid of $\q(4,q)$. Using $\q(4,q) \cong T_2(\C)$ if and only if
$\C$ is a conic, $\oo$ is equivalent with a set of points $U \cup \{\infty\}$, $U$ consisting
exclusively of points of type (i), since any point of $\q(4,q)$ can play the role of the point
$\infty$. Since no two points of $\oo$ are collinear in $\q(4,q)$, no two points of $U$ may
determine a line of type (a) of $T_2(\C)$, and since all lines of type (a) meet $\pi$ in a point
of $\C$, a (partial) ovoid of $\q(4,q)$ of size $q^2+1-\epsilon$ is equivalent
to a set $U$ of $q^2-\epsilon$ points of $\ag(3,q)$, not determining the points of a conic
at infinity. 

An immediate application of Theorem~\ref{th:ball2008} is the following theorem.

\begin{theorem}[\cite{Ball:2008rc}, Section 4]
Let $U$ be a set of $q^2$ points in $\ag(3,p)$, $p$ prime, whose non determined directions contain a conic.
Then every plane of $\ag(3,q)$ meets $U$ in $0 \mod p$ points.
\end{theorem}

And the previous theorem has the classification of all ovoids of $\q(4,p)$, $p$ prime,
as a consequence.

\begin{cor}[\cite{Ball:2008rc}, Section 4]
An ovoid of $\q(4,p)$, $p$ prime, is necessarily contained in a hyperplane of $\pg(4,q)$, meeting
$\q(4,q)$ in the points of an elliptic quadric.
\end{cor}

Note that the previous corollary and theorem were first proved in \cite{Ball:2004vn}, using a purely algebraic approach that
is unrelated to direction problems.

The question addressed in \cite{De-Beule:2008fk} is whether a partial ovoid of $\q(4,q)$
of size $q^2-1$ can be maximal, i.e. whether the corresponding set $U$ of size $q^2-2$ 
can be non-extendable. As an immediate application of Theorem~\ref{th:dbgacs}, the following theorem 
is described in \cite{De-Beule:2008fk}.

\begin{theorem}[\cite{De-Beule:2008fk}, Theorem 3]
Let $q=p^h$, $p$ odd prime, $h > 1$. Then $\q(4,q)$ has no maximal partial ovoids of size $q^2-1$.
\end{theorem}

A (finite) \emph{partial geometry} is a point-line geometry $\IS$ that is a generalization of a GQ. 
To define a partial geometry, the third GQ axiom is replaced by the following:
\begin{compactenum}[(iii)]
\item There exists a fixed integer $\alpha > 0$, such that if $x$ is a point and 
$L$ is a line not incident with $x$, then there are exactly $\alpha$ pairs 
$(y_i,M_i) \in \PP \times \B$ for which $x \I M_i \I y_i \I L$.
\end{compactenum}
The integers $s$, $t$ and $\alpha$ are the parameters of $\cS$. The \emph{dual} $\cS^D$ of a
partial geometry $\IS$ is the incidence structure $(\B,\PP,\I)$. It is a partial geometry
with parameters $s^D = t$, $t^D = s$, $\alpha^D = \alpha$. Clearly, a partial geometry
with $\alpha=1$ is a GQ.

We need some special pointsets in $\pg(2,q)$ to describe our favorite partial geometries. 
An {\em arc of degree $d$} of $\pg(2,q)$ is a set $\K$ of points such that
every line of $\pg(2,q)$ meets $\K$ in at most $d$ points. If $\K$ contains $k$ points,
that it can also be called a $\{k,d\}$-arc. A typical example is a conic,
which is a $\{q+1,2\}$-arc. The size of an arc of degree $d$ can not exceed $dq-q+d$.
A $\{k,d\}$-arc $\K$ for which $k=dq-q+d$, or equivalently, such that every line that meets $\K$,
meets $\K$ in exactly $d$ points, is called {\em maximal}. With this definition, a conic
is a non maximal $\{q+1,2\}$-arc, and it is well known that if $q$ is even, a conic, together with its 
nucleus, is a $\{q+2,2\}$-arc, which is complete. We mention that a $\{q+1,2\}$-arc is also called
an {\em oval}, and a $\{q+2,2\}$-arc is also called a {\em hyperoval}. When $q$ is odd, all ovals are
conics, and no $\{q+2,2\}$-arcs exist. Let $q=2^h$, then every oval has a nucleus, and 
so can be extended to a hyperoval. Much more examples of hyperovals, different from a conic and its 
nucleus, are known. Maximal $\{k,d\}$-arcs exist for $d=2^e$, $1 \leq e \leq h$. Several infinite families and 
constructions are known. We refer to \cite{DCD:2012} for an overview, and detailed references to
mentioned results here.

Let $q$ be even and let $\K$ be a maximal $\{k,d\}$-arc of $\pg(2,q)$. We define the incidence structure $T_2^*(\K)$ as
follows. Embed $\pg(2,q)$ as a hyperplane $H_{\infty}$ in $\pg(3,q)$. The points of $\cS$ are the points of $\pg(3,q)\setminus H_{\infty}$. The 
lines of $\cS$ are the lines of $\pg(3,q)$ not contained in $H_{\infty}$, and meeting $H_{\infty}$ in a point of $\K$. The incidence
is the natural incidence of $\pg(3,q)$. One can check easily, using that $\K$ is a maximal $\{k,d\}$-arc, that $T_2^*(\K)$
is a partial geometry with parameters $s=q-1$, $t=k-1=(d-1)(q+1)$, and $\alpha = d-1$.

The definitions of (maximal) (partial) ovoid can be taken over for partial geometries. As in the GQ case, a (maximal) (partial)
ovoid of $T_2(\K)$ is equivalent with a (extendable) set of points $U$ of $\ag(3,q)$ not determining the points of $\K$ at infinity.
A maximal $\{k,d\}$-arc can not be contained in a conic. Therefore Theorem~\ref{th:dbszt2013} yields almost immediately 
the following theorem as a corollary.

\begin{theorem}[\cite{De-Beule:2013uq}, Corollary 18]
Let $\B$ be a partial ovoid of size $q^2-2$ of the partial geometry $T_2^*(\K)$, then $\B$ is always extendable to an ovoid. 
\end{theorem}

Note that the following theorem is found in \cite{Ball:2006mz}. It is obtained indirectly as a corollary of Theorem~\ref{th:ag3q1}
and some extra work.

\begin{theorem}[\cite{Ball:2006mz}, Corollary 3.2]
Let $\cH$ be a hyperoval of $\pg(2,q)$, $q$ even. Let $\oo$ be an ovoid of $T_2^*(\cH)$, then every plane of $\pg(3,q)$ meets
$\oo$ in an even number of points. Moreover two planes meeting $\pg(3,q) \setminus \ag(3,q)$ in the same line intersect
$\oo$ either both in $0 \mod 4$ points or both in $2 \mod 4$ points.
\end{theorem}

\section*{Acknowledgement}

The author thanks Tam\'as Sz\H{o}nyi for many interesting remarks on the first version of this paper.


\noindent Jan De Beule\\
Ghent University\\
Department of Mathematics\\
Krijgslaan 281, S22\\
B-9000 Gent\\
Belgium\\
\\
and\\
\\
Vrije Universiteit Brussel\\
Department of Mathematics\\
Pleinlaan 2\\
B-1050 Brussel\\
Belgium\\
\\
\url{jdebeule@cage.ugent.be}
\end{document}